\documentclass[12pt]{amsproc}
\usepackage{amsmath,amssymb,amscd,stmaryrd}
\usepackage[mathscr]{euscript}

\newcommand{\sA}{\mathscr A}
\newcommand{\sB}{\mathscr B}
\newcommand{\sC}{\mathscr C}
\newcommand{\sD}{\mathscr D}
\newcommand{\sE}{\mathscr E}
\newcommand{\sF}{\mathscr F}
\newcommand{\sG}{\mathscr G}
\newcommand{\sH}{\mathscr H}
\newcommand{\sI}{\mathscr I}
\newcommand{\sK}{\mathscr K}
\newcommand{\sM}{\mathscr M}
\newcommand{\sN}{\mathscr N}
\newcommand{\sP}{\mathscr P}
\newcommand{\sR}{\mathscr R}
\newcommand{\sU}{\mathscr U}
\newcommand{\sX}{\mathscr X}

\begin{document}

\title[Commutants mod Normed Ideals]{Commutants mod Normed Ideals}
\author{Dan-Virgil Voiculescu}
\address{D.V. Voiculescu \\ Department of Mathematics \\ University of California at Berkeley \\ Berkeley, CA\ \ 94720-3840}
\thanks{Supported in part by NSF Grant DMS-1665534. Part of this paper was written while visiting IPAM in Spring 2018 for the Quantitative Linear Algebra program, which was supported by a NSF grant.}
\dedicatory{Dedicated to Alain Connes on the occasion of his 70th birthday.}


\maketitle

\section{Introduction}
\label{sec1}

The full name of the normed ideals to which we refer in the title, is ``symmetrically normed ideals of compact operators'', among which the Schatten--von~Neumann $p$-classes are the most familiar. The commutants are commutants modulo a normed ideal of $n$-tuples of self-adjoint operators, or equivalently of the algebras generated by these operators. When the normed ideal is $\sK$, that is the ideal of compact operators, then up to factoring by $\sK$, this is roughly how one arrives at the Paschke dual of a finitely generated $C^*$-algebra (\cite{35}, see also \cite{25}), that is a basic duality construction in the $K$-theory of $C^*$-algebras, that in essence hails from the abstract elliptic operators of Atiyah \cite{2}. What happens if $\sK$ is replaced by a smaller normed ideal $\sI$? As we will see, commutants mod~$\sI$ are not simply ``smooth versions'' of those mod~$\sK$, and they are often closer to $C^*$-algebras than one would expect, while they connect with many questions in perturbation theory.

From an operator theory perspective, following the Brown--Douglas--Fillmore work \cite{7}, \cite{8} and the development of the $K$-theory of $C^*$-algebras all the way to Kasparov's bivariant theory (\cite{26}, see also \cite{5}, \cite{23}) many things about compact perturbations can now be understood from this point of view.

When other normed ideals than $\sK$ are considered, there are different aspects to be taken into account. For instance, Alain Connes' cyclic cohomology is the algebraic framework for trace-formulas like the Helton--Howe formula for almost normal operators \cite{22}, \cite{15}. On the other hand, the invariance of Lebesgue absolutely continuous spectra under trace-class perturbations, fits with our modulus of quasicentral approximation invariant which is one focus of the present article.

Starting with an adaptation \cite{43} of the Voiculescu non-commutative Weyl--von~Neumann type theorem \cite{42} to normed ideals other than $\sK$, we found that $k_{\sI}(\tau)$, the modulus of quasicentral approximation for the $n$-tuple $\tau$ relative to the normed ideal $\sI$, underlies many questions concerning perturbation of operators \cite{43}, \cite{44}, \cite{49}. This quantity has also turned out to be connected to the Kolmogorov--Sinai dynamical entropy ad to the supramenability of groups \cite{47}, \cite{54}. We should also mention that, as we found in \cite{43}, \cite{44}, \cite{45}, in many of these questions Lorentz $(p,1)$ ideals instead of Schatten--von~Neumann $p$-classes give sharp results when $p > 1$.

To Alain Connes' non-commutative geometry, the normed ideals are purveyors of infinitesimals \cite{15}. The machinery around $k_{\sI}(\tau)$ has found technical uses in the spectral characterization of compact manifolds \cite{16} (i.e., their characterization as non-commutative manifolds) and in results about unbounded Fredholm modules \cite{13}.

Since the numerical invariant $k_{\sI}(\tau)$ turned out to have good properties and to unify several perturbation problems, we recently looked whether there is not actually more structure around this number. The commutant $\sE(\tau;\sI)$ of an $n$-tuple of operators $\tau$~mod a normed ideal $\sI$ provides such structure. The first step in the study of $\sE(\tau;\sI)$ is to introduce $\sK(\tau;\sI)$ its ideal of compact operators and the corresponding Calkin algebra $\sE/\sK(\tau;\sI)$. Under suitable assumptions on $k_{\sI}(\tau)$, there are many functional analysis similarities between $(\sK(\tau;\sI),\sE(\tau;\sI),\sE/\sK(\tau;\sI))$ and the usual $(\sK,\sB,\sB/\sK)$ where $\sB$ denotes the bounded operators. With this analogy as a guide we also took the first steps in computing some $K$-groups of such algebras. The $K$-theory can be quite rich and this demonstrates why the $\sE/\sK(\tau;\sI)$ should not be thought as being ``smooth versions'' of $\sE/\sK(\tau;\sK)$, i.e., essentially of Paschke duals.

As a general comment about algebras associated to perturbations from a normed ideal $\sI$, it appears that $\sE(\tau;\sI)$, which is closer to the abstract elliptic operators point of view of Atiyah, being a Banach algebra may have from a functional analysis point of view some advantages over dealing with homomorphisms of algebras into $B/\sI$ in the Brown--Douglas--Fillmore style. The algebra $B/\sI$ is not a Banach algebra, though smooth functional calculus of various kinds can still be performed in $B/\sI$. On the other hand $\sE/\sK(\tau;\sI)$ if $k_{\sI}(\tau) < \infty$ is actually isomorphic to a $C^*$-algebra.

Very recently we found that $k_{\sI}(\tau)$ remains an effective tool also when generalized to handling of hybrid perturbations. That is instead of a normed ideal $\sI$, we will have an $n$-tuple of normed ideals $(\sI_1,\dots,\sI_n)$ and the perturbation of $\tau$ to $\tau'_1$ is such that $T_j - T'_j \in \sI-j$, $1 \le j \le n$. We will also very briefly mention a few results in this direction \cite{57}, \cite{58}, which are quite sharp.

Here is in brief how the survey proceeds. After this introduction we provide some background on the normed ideals we will use. We then as a motivating example recall the main facts about normed ideal perturbations of one self-adjoint operator. Then we introduce the invariant $k_{\sI}(\tau)$ and its basic properties. Next we give a version of the Voiculescu theorem adapted to normed ideals. After this we pass to the applications of this machinery to normed ideal perturbations of $n$-tuples of commuting hermitian operators. Next we explain the endpoint properties of $k^-_{\infty}(\tau)$ which is the case when $\sI$ is the Macaev ideal. The way the Kolmogorov--Sinai entropy is related to $k^-_{\infty}(\tau)$ is explained after this. Then we discuss results for finitely generated groups and the result and open problem about $k^-_{\infty}$ and supramenable groups. After this, we go over to commutants mod normed ideals. We explain that $k_{\sI}(\tau)$ is related to approximate units of the compact ideal $\sK(\tau;\sI)$ of $\sE(\tau;\sI)$. The Banach space properties of $\sE(\tau;\sI)$ are then discussed as well as properties of $\sE/\sK(\tau;\sI)$ where $k_{\sI}(\tau) < \infty$ or $k_{\sI}(\tau) = 0$. We then look at the results about $K_0(\sE(\tau;\sI))$ for $n$-tuples of commuting hermitian operators, which, to simplify matters we choose here to be multiplication operators by the coordinate functions in the $L^2$-space with respect to Lebesgue measure on a hypercube. We then discuss a few results in the hybrid setting. Quite briefly finally some applications to unbounded Fredholm modules are then pointed out.

The study of commutants mod normed ideals and of the invariant $k_{\sI}(\tau)$ and of the hybrid generalization are still at an early stage and from a reading of this survey one realizes the multitude of open problems of varying degrees of difficulty. Still at the end we briefly mention some sample open problems which had not appeared with the presentation of results.

At various points in this exposition, in order to avoid technical detail, we did not aim at the most general or technically strongest version of the results. We hope this kind of simplification will make it easier for the reader to focus on the big picture.

\section{Background on normed ideals}
\label{sec2}

Throughout we shall denote by $\sH$ an infinite dimensional separable complex Hilbert space and by $\sB(\sH)$, $\sK(\sH)$, $\sR(\sH)$ or simply $\sB,\sK,\sR$ the bounded operators, the compact operators and the finite rank operators on $\sH$. The Calkin algebra is then $\sB/\sK$ and $p: \sB \to \sB/\sK$ will be the canonical homomorphism. A normed ideal $(\sI,|\ |_{\sI})$ is an ideal of $\sB$ so that $\sR \subset \sI \subset \sK$ and $|\ |_{\sI}$ is a Banach space norm on $\sI$ satisfying a set of axioms which can be found in \cite{20}, \cite{40}. In particular if $A,B \in \sB$ and $X \in \sI$ we have $|AXB|_{\sI} \le \|A\|\,|X|_{\sI}\,\|B\|$. If $\sI$ is a normed ideal, its dual as a Banach space is also a normed ideal $\sI^d$, where now we have to also allow the possibility that $\sI^d = \sB$, and the duality is given by $(X,Y) \to \mbox{Tr}\ XY$. If $\sI$ is a normed ideal then the closure of $\sR$ in $\sI$ is also a normed ideal which we shall denote by $\sI^{(0)}$ and which may be strictly smaller than $\sI$.

If $1 \le p < \infty$, the Schatten--von~Neumann $p$-class $\sC_p$ is the normed ideal of operators $X \in \sK$ so that $|X|_p = (\mbox{Tr}((X^*X)^{p/2}))^{1/p} < \infty$ endowed with the norm $|\ |_p$. In particular, $\sC_1$ is the trace-class and $\sC_2$ is the ideal of Hilbert--Schmidt operator. Another scale of normed-ideals which we shall use here is $(\sC^-_p,|\ |^-_p)$, $1 \le p \le \infty$ which are the Lorentz $(p,1)$-ideals. If $s_1 \ge s_2 \ge \dots$ are the eigenvalues of the compact operator $(X^*X)^{1/2}$ then $|X|^-_p = \sum_{j \in {\mathbb N}} s_jj^{-1+1/p}$. Of particular interest is the case where $p = \infty$ and $|X|^-_{\infty} = \sum_{j \in {\mathbb N}} s_jj^{-1}$. If $p = 1$, $\sC^-_1 = \sC_1$, but otherwise $\sC_p^- \subset \sC_p$ but $\sC_p^- \ne \sC_p$, $1 < p \le \infty$, where $\sC_{\infty} = \sK$. The ideal $\sC^-_p$ can also be described as the smallest normed ideal for which the norm on projections is equivalent to the $p$-norm, when $1 \le p < \infty$ while $|P|^-_{\infty} \sim \log \mbox{Tr } P$ when $p = \infty$ for a projection $P$. The ideal $\sC^-_{\infty}$ is also called the Macaev ideal.

The normed ideals can also be viewed as the non-commutative analogue of classical Banach sequence spaces \cite{29}.

\section{The theorems of Weyl--von~Neumann--Kuroda and of Kato--Rosenblum}
\label{sec3}

Let $A$ and $B$ be hermitian operators on $\sH$ and let $(\sI,|\ |_{\sI})$ be a normed ideal. Recall that $p$ denotes the homomorphism onto the Calkin algebra. Thus the essential spectrum $\sigma(p(A))$ of the hermitian operator $A$ is obtained from its spectrum $\sigma(A)$ by removing the isolated points $\lambda \in \sigma(A)$ which corresponds to eigenvalues of finite multiplicity. By results of Weyl, von~Neumann and Kuroda, see \cite{27}, we have:

\bigskip
\begin{center}
\begin{minipage}{4in}
Assume $\sI \ne \sC_1$, $\varepsilon > 0$, and $\sigma(A) = \sigma(B) = \sigma(p(A)) = \sigma(p(B))$. Then there is a unitary operator $U$ so that
\[
|UBU^* - A|_{\sI} < \varepsilon.
\]
\end{minipage}
\end{center}

\bigskip
Note that, given $A$, we may choose $B$ to be diagonal in an orthonormal basis and thus we get $|X-A|_{\sI} < \varepsilon$ where $X = UBU^*$ is an operator which can be diagonalized in an orthonormal basis.

The trace-class $\sC_1$ is actually the smallest normed ideal. If $\sI = \sC_1$, the previous result fails because of the Lebesgue absolutely continuous spectrum, which is a conserved quantity under trace-class perturbations. This is a consequence of the Kato--Rosenblum theorem of abstract scattering theory (\cite{27}, \cite{37}). If $X = X^*$, we say that $X$ has Lebesgue absolutely continuous spectrum, if its spectral measure $E(X;\cdot)$ is absolutely continuous w\.r\.t\. Lebesgue measure (equivalently this is that the scalar measures $\langle E(X;\cdot)\xi,\eta\rangle$ for all $\xi,\eta \in \sH$ are Lebesgue absolutely continuous). Given $X$, the Hilbert space $\sH$ splits in a unique way $\sH = \sH_{ac} \oplus \sH_{sing}$ into $X$-invariant subspaces so that $X\mid \sH_{ac}$ has Lebesgue absolutely continuous spectrum, while $X \mid \sH_{sing}$ has singular spectrum, that is the spectral measure of $X \mid \sH_{sing}$ is carried by a Borel set of Lebesgue measure zero. It is a corollary of the Kato--Rosenblum theorem that:

\bigskip
\begin{center}
\begin{minipage}{4in}
if $X - A \in \sC_1$ then $X \mid \sH_{ac}(X)$ and $A \mid \sH_{ac}(A)$ are unitarily equivalent.
\end{minipage}
\end{center}

\bigskip
The Kato--Rosenblum theorem actually provides two intertwiners, to achieve the unitary equivalence, the generalized wave operators $W_{\pm}$:

\bigskip
\begin{center}
\begin{minipage}{4in}
under the assumption $X - A \in \sC_1$ the strong limits
\[
W_{\pm} = s - \lim_{t \to \pm \infty} e^{itA}e^{-itX} \mid \sH_{ac}(X)
\]
exist. Moreover we have $W_{\pm}\sH_{ac}(X) = \sH_{ac}(A)$ and
\[
W_{\pm}(X \mid \sH_{ac}(X)) = (A \mid \sH_{ac}(A))W_{\pm}.
\]
\end{minipage}
\end{center}

\bigskip
Thus the Lebesgue absolutely continuous part $A \mid \sH_{ac}(A)$ is conserved up to unitary equivalence under perturbations in $\sC_1$ and cannot be diagonalized.

On the other hand for the singular spectrum, the ideal $\sC_1$ is not different from other ideals. For instance the Weyl--von~Neumann--Kuroda theorem holds also for $\sC_1$ if the spectrum is singular that is if $\sH_{ac}(A) = \sH_{ac}(B) = 0$.

Concerning what goes into the proofs of these theorems, the Weyl--von~Neumann--Kuroda results rely essentially on partitioning
\[
I = E(A;\omega_1) + \dots + E(A;\omega_n)
\]
where $\omega_1,\dots,\omega_n$ is a partition of $\sigma(A)$ into Borel sets of small diameter, while the Kato--Rosenblum theorem uses some Fourier analysis, which can be viewed as related to the $L^2$-boundedness of the Hilbert transform.

To conclude this discussion we should also say that the role of assumptions about spectra and essential spectra $\sigma(A) = \sigma(p(A))$ will become clearer in the next section when we consider $C^*$-algebras.

\section{The theorem of Voiculescu}
\label{sec4}

After the work of Brown--Douglas--Fillmore which solved the unitary conjugacy question for normal elements of the Calkin algebra, it became clear that it is preferable in this kind of question to view operators or $n$-tuples of operators as representations of the $C^*$-algebra which they generate. For instance given an $n$-tuple of commuting hermitian operators $A_1,\dots,A_n$ this is just the representation of the commutative $C^*$-algebra $C(K)$ of continuous functions on their joint spectrum $K = \sigma(A_1,\dots,A_n)$ which arises from functional calculus $\rho(f) = f(A_1,\dots,A_n)$.

One version of the Voiculescu theorem (\cite{42}, see also \cite{1}) is a non-commutative generalization of the perturbation results for one hermitian operator in case $\sI = \sK$, the ideal of compact operators:

\bigskip
\begin{center}
\begin{minipage}{4in}
If $\sA$ is a unital separable $C^*$-algebra and $\rho_1,\rho_2: \sA \to \sB(\sH)$ are unital $*$-homomorphisms such that $\ker p \circ \rho_1 = \ker p \circ \rho_2 = 0$, then there is a unitary operator $U$ so that $\rho_1(a) - U\rho_2(a)U^*(a) \in \sK$ for all $a \in \sA$.
\end{minipage}
\end{center}

\bigskip
For instance to recover the Weyl--von~Neumann result for two hermitian operators $A,B$ with $\sigma(A) = \sigma(B) = \sigma(p(A)) = \sigma(p(B)) = K$ one takes $\sA = C(K)$, $\rho_1(f) = f(A)$, $\rho_2(f) = f(B)$ and the conclusion is applied to the particular choice of $f$ being the identical function.

This very general result for compact perturbations, obviously leads to the question: What happens when the compact operators are replaced by a smaller normed ideal $\sI$? We found that there is a key quantity which needs to be taken into account and which will be discussed in the next section.

\section{The invariant $k_{\sI}(\tau)$}
\label{sec5}

Let $\tau = (T_j)_{1 \le j \le n}$ be an $n$-tuple of bounded operators on $\sH$ and $(\sI,|\ |_{\sI})$ a normed ideal. Then the {\em modulus of quasicentral approximation} is defined as follows (\cite{43}, \cite{45}):

\bigskip
\begin{center}
\begin{minipage}{4in}
$k_{\sI}(\tau)$ is the least $C \in [0,\infty]$, such that there exist finite rank operators $0 \le A_m \le I$ so that $A_m \uparrow I$ and we have
\[
\lim_{m \to \infty} \max_{1 \le j \le n} |[A_m,T_j]|_{\sI} = C.
\]
\end{minipage}
\end{center}

\bigskip
If $\sI = \sC_p$ or $\sI = \sC_p^-$ we denote $k_{\sI}(\tau)$ by $k_p(\tau)$ or respectively by $k^-_p(\tau)$. In particular in the case of the Macaev ideal $\sI = \sC^-_{\infty}$ we get $k_{\infty}^-(\tau)$. Here are some of the first properties of this invariant (see \cite{45}).

\bigskip
\begin{center}
\begin{minipage}{4in}
$1^{\circ}$ $[1,\infty] \ni p \to k_p^-(\tau) \in [0,\infty]$ is a decreasing function of $p$.
\end{minipage}

\medskip\smallskip
\begin{minipage}{4in}
$2^{\circ}$ for a given $\tau$ there is $p_0 \in [1,\infty]$ so that if $p \in [1,p_0)$ then $k^-_p(\tau) = \infty$, while if $p \in (p_0,\infty]$ then $k^-_p(\tau) = 0$. (Note that if $k^-_p(\tau) \in (0,\infty)$ then we must have $p_0 = p$.)
\end{minipage}

\medskip\smallskip
\begin{minipage}{4in}
$3^{\circ}$ assuming $\sR$ is dense in $\sI$, then if $\tau,\tau'$ are $n$-tuples so that $T_j-T'_j \in \sI$, $1 \le j \le n$, then $k_{\sI}(\tau) = k_{\sI}(\tau')$.
\end{minipage}

\medskip\smallskip
\begin{minipage}{4in}
$4^{\circ}$ assuming $\tau = \tau^*$ and that $\sR$ is dense in $\sI$, we have: $k_{\sI}(\tau) > 0$ iff there exist $Y_j = Y^*_j \in \sI^{\mbox{dual}}$, $1 \le j \le n$ so that $i\sum_j [T_j,Y_j] \in \sC_1 + \sB(\sH)_+$ and $\mbox{Tr}\ i \sum_j [T_j,Y_j] > 0$.
\end{minipage}
\end{center}

\bigskip
To prove that $k_{\sI}(\tau) = 0$ or that $k_{\sI}(\tau) < \infty$ one can use the definition of $k_{\sI}(\tau)$ and find a suitable sequence of operators $A_n$. To prove that $k_{\sI}(\tau) > 0$ is usually more difficult as one has to find suitable $Y_j$ $(1 \le j \le n)$ which satisfy property $4^{\circ}$. For instance in the case of $\sI = \sC_1$ and $n = 1$, this boils down to the boundedness of the Hilbert transform in $L^2$. Indeed, remark that if $T$ is the multiplication operator in $L^2([0,1],d\lambda)$ by the coordinate function and $H$ is the compression of the Hilbert transform to $L^2([0,1],d\lambda)$ then $[T,H] = iP$ where $P$ is the rank one projection onto the constant functions.

The reader may have wondered why we did not pay more attention to $k_p(\tau)$ instead of focusing on $k^-_p(\tau)$. The reason is that if $p = 1$ then $\sC^-_1 = \sC_1$ and we have $k_1(\tau) = k^-_1(\tau)$, while otherwise, if $p > 1$ we have that $k_p(\tau) \in \{0,\infty\}$ (\cite{45}).

\section{Some uses of $k_{\sI}(\tau)$}
\label{sec6}

A general result that uses $k_{\sI}(\tau)$ is the adaptation of the Voiculescu theorem \cite{43}, to deal with other normed ideals than $\sK$. Here is a version of such an adaptation:

\bigskip
\begin{center}
\begin{minipage}{4in}
Let $\sA$ be a $C^*$-algebra generated by $X_k$, $1 \le k \le n$ and let $\rho_1,\rho_2$ be unital $*$-representations on $\sH$ of $\sA$, so that $\ker p \circ \rho_j = 0$, $j = 1,2$. Assume moreover that $k_{\sI}(\rho_j(X_k)_{1 \le k \le n}) = 0$, $j = 1,2$. Then there is a unitary operator $U$ so that
\[
U\rho_1(X_k)U^* - \rho_2(X_k) \in \sI,\ 1 \le k \le n.
\]
\end{minipage}
\end{center}

\bigskip
Actually, this version which we chose for simplicity, is a corollary of a more general absorption version (see \cite{43}). Note also that one can also have that
\[
|U\rho_1(X_k)U^* - \rho_2(X_k)|_{\sI} < \varepsilon,\ 1 \le k \le n
\]
for a given $\varepsilon > 0$.

Remark that the previous result when applied to $\sA = C(K)$ where $K \subset {\mathbb R}$ is a compact set and $X$ is the identical function $X(t) = t$, reduces the Weyl--von~Neumann--Kuroda theorem to proving the following fact:

\bigskip
\begin{center}
\begin{minipage}{4in}
if $T = T^*$ and $\sI \ne \sC_1$, then $k_{\sI}(T) = 0$.
\end{minipage}
\end{center}

\bigskip
There is also the following general construction for a given $n$-tuple $\tau$ and a normed ideal $\sI$.

If $(P_i)_{i \in I}$ are projections in the commutant $(\tau)'$ so that $k_{\sI}(\tau \mid P_i\sH) = 0$, then also $P = \underset{i \in I}{\bigvee} P_i$ is so that $k_{\sI}(\tau \mid P\sH) = 0$. In particular, there is a largest reducing subspace $\sH_s$ of $\tau$, so that $k_{\sI}(\tau \mid \sH_s) = 0$. The subspace $\sH_s$ is called the $\sI$-singular subspace of $\tau$, while $\sH_s = \sH \ominus \sH_s$ is called the $\sI$-absolutely continuous subspace of $\tau$.

The names given to these reducing subspaces of $\tau$ are motivated by the case of commuting hermitian operators. In particular, in the simplest case of one hermitian operator $T = T^*$ and $\sI = \sC_1$, the $\sC_1$-singular and $\sC_1$-absolutely continuous subspaces of $T$ are precisely the Lebesgue singular and Lebesgue absolutely continuous subspaces of $T$. The use of $k_1$ reduces this to the following fact:

\bigskip
\begin{center}
\begin{minipage}{4in}
$k_1(T) = 0 \Leftrightarrow$ the spectral measure of $T$ is singular w.r.t. Lebesgue measure.
\end{minipage}
\end{center}

\bigskip
Note that $\Rightarrow$ relies essentially on the $L^2$-boundedness of the Hilbert transform.

Much of our initial motivation for developing a machinery based on the invariant $k_{\sI}(\tau)$ for studying perturbations of $n$-tuples of operators was to extend perturbation results for one hermitian operator to commuting $n$-tuples of hermitian operators. For instance for $n = 2$, which is equivalent to dealing with one normal operator $N$, there was the problem attributed to P.~R. Halmos whether $N = \sD + K$ where $\sD$ was diagonalizable and $K \in \sC_2$. In essence, this problem was whether $\sC_2$ plays the same role for pairs of commuting hermitian operators that $\sC_1$ plays for singletons. It turned out that normal operators can be diagonalized $\mod \sC_2$ (\cite{43}), but it is better to go over to $n$-tuples of commuting hermitian operators and describe the general results obtained with the modulus of quasicentral approximation.

\section{Perturbations of commuting $n$-tuples of hermitian operators}
\label{sec7}

It turned out that for $n$-tuples of commuting hermitian operators there is a threshold ideal $\sC^-_n$. This also explains the $\mod \sC_2$ diagonalization of normal operators. Here are some of the main results for such $n$-tuples (\cite{3}, \cite{7}, \cite{43}, \cite{44}).

\bigskip
\begin{center}
\begin{minipage}{4in}
$1^{\circ}$ If $\sI$ is a normed ideal and $\sI \supset \sC^-_n$, $\sI \ne \sC^-_n$ and $\tau$ is an $n$-tuple of commuting hermitian operators then $k_{\sI}(\tau) = 0$. In particular, there is a diagonalizable $n$-tuple $\delta$ so that $\tau \equiv \delta \mod \sI$.
\end{minipage}

\medskip\smallskip
\begin{minipage}{4in}
$2^{\circ}$ If $\tau$ and $\tau'$ are $n$-tuples of commuting hermitian operators and $\tau \equiv \tau' \mod \sC^-_n$, then their Lebesgue absolutely continuous parts $\tau_{ac}$ and $\tau'_{ac}$ are unitarily equivalent.
\end{minipage}

\medskip\smallskip
\begin{minipage}{4in}
$3^{\circ}$ There is a universal constant $0 < \gamma_n < \infty$ so that if $\tau$ is an $n$-tuple of commuting hermitian operators, then
\[
(k^-_n(\tau))^n = \gamma_n \int_{{\mathbb R}^n} m(s)d\lambda(s)
\]
where $\lambda$ is Lebesgue measure and $m(s)$ is the multiplicity function of the Lebesgue absolutely continuous part of $\tau$. If $n = 1$, then $\gamma_1 = 1/\pi$.
\end{minipage}
\end{center}

\bigskip
The Kato--Rosenblum theorem was also generalized using $k^-_n(\tau)$, but the results (see \cite{44}) are perhaps not complete. For $n \ge 3$ we have a very general result showing that some very general generalized wave operators exist and are actually all equal. For $n = 2$ we get the existence of enough non-trivial intertwiners, but the convergence is not strong. Note that we don't have a proof of the usual Kato--Rosenblum theorem using $k_1(\tau)$. Its corollary about the unitary equivalence of absolutely continuous parts can be however recovered in case the multiplicity function of one of the operators is integrable by using the formula for $k_1$ applied to $f(A)$ for $C^{\infty}$-functions $f$.

It is an open problem whether the very strong Kato--Rosenblum type results we proved for $n \ge 3$ in \cite{44} also hold for $n = 2$.

The machinery based on $k_{\sI}(\tau)$ for dealing with perturbations of commuting $n$-tuples of hermitian operators applies as soon as we know the decomposition $\sH = \sH_a(\tau) \oplus \sH_s(\tau)$ for a given normed ideal $\sI$. This essentially means to study the $n$-tuple of multiplication operators by the coordinate functions in $L^2({\mathbb R}^n,\mu)$ where $\mu$ is a more general measure, for instance more like a Hausdorff $p$-dimensional measure $1 \le p < n$, $p$ not necessarily an integer. The essential difficulty is in showing that certain singular integrals give operators in certain normed ideals in order to show that $k^-_p(\tau) > 0$. Rather general results of this kind were obtained in our joint work with Guy David \cite{17}. Here is the key result from \cite{17}:

\bigskip
\begin{center}
\begin{minipage}{4in}
{\em Let $\mu$ be a Radon probability measure with compact support on ${\mathbb R}^n$ so that the Ahlfors condition}
\[
\mu(B(x,r)) \le Cr^p,\ \forall x \in {\mathbb R}^n,\ r \le 1
\]
{\em holds for a certain $p > 1$. Let further $\tau_{\mu}$ be the $n$-tuple of multiplication operators by the coordinate functions in $L^2({\mathbb R}^n,d\mu)$. Then we have:}
\[
k^-_p(\tau_{\mu}) > 0.
\]
\end{minipage}
\end{center}

\section{$k^-_p(\tau)$ at the endpoint $p = \infty$ and dynamical entropy}
\label{sec8}

In case $\sI = \sC^-_{\infty}$, the Macaev ideal, the invariant $k^-_{\infty}(\tau)$ has remarkable properties (\cite{45}).

\bigskip
\begin{center}
\begin{minipage}{4in}
{\em Let $\tau$ be an $n$-tuple of bounded operators. We have:}
\end{minipage}

\medskip\smallskip
\begin{minipage}{4in}
(i) {\em $k_{\infty}^-(\tau) < \infty$, more precisely $k_{\infty}^-(\tau) \le 2\|\tau\| \log (2n+1)$}
\end{minipage}

\medskip\smallskip
\begin{minipage}{4in}
(ii) $k_{\infty}^-(\tau) = k_{\infty}^-(\tau \otimes I_{\sH})$
\end{minipage}

\medskip\smallskip
\begin{minipage}{4in}
(iii) {\em if $\sI$ is a normed ideal so that $\sI \supset \sC^-_{\infty}$ and $\sI \ne \sC^-_{\infty}$ then}
\[
k_{\sI}(\tau) = 0
\]
\end{minipage}

\medskip\smallskip
\begin{minipage}{4in}
(iv) {\em if $S_1,\dots,S_n$ are isometries with orthogonal ranges and $n \ge 2$ then}
\[
k^-_{\infty}(S_1,\dots,S_n) > 0.
\]
\end{minipage}
\end{center}

\bigskip
We saw that $k^-_p$ for finite integer $p$ is related to $p$-dimensional Lebesgue measure and somewhat more loosely when $p$ is not an integer to corresponding quantities of Hausdorff dimension $p$. When $p = \infty$ we have found that instead of a $p$-dimensional measure there are connections to dynamical entropy. Using $k^-_{\infty}$ a quantity ``approximately'' equivalent to the Kolmogorov--Sinai dynamical entropy can be obtained (\cite{47}). Here is how this {\em dynamical perturbation entropy} (\cite{47}) is constructed.

Let $\theta$ be a measure-preserving automorphism of a probability measure space $(\Omega,\Sigma,\mu)$, $\mu(\Omega) = 1$. Let further $U_{\theta}$ be the unitary operator in $L^2(\Omega,\Sigma,\mu)$ induced by $\theta$ and $\Phi$ the set of multiplication operators in $L^2(\Omega,\Sigma,\mu)$ by measurable numerical functions which take finitely many values. The dynamical perturbation entropy is defined by the formula

\bigskip
\begin{center}
\begin{minipage}{4in}
\[
\sH_p(\theta) = \sup_{\begin{matrix}
\varphi \subset \Phi \\
\varphi \text{ finite}
\end{matrix}} k^-_{\infty}(\varphi \cup \{U_{\theta}\}).
\]
\end{minipage}
\end{center}

\bigskip
\noindent
This is the definition from \cite{49}. It is easy to show that it is equal to the quantity denoted by ${\tilde \sH}_{P}(\theta)$ in \cite{47}, \cite{48}.

Comparing $\sH_{P}(\theta)$ to the Kolmogorov--Sinai entropy $h(\theta)$ we have the following results \cite{48}.

\bigskip
\begin{center}
\begin{minipage}{4in}
(i) There are universal constants $0 < C_1 < C_2 < \infty$ so that
\[
C_1h(\theta) \le \sH_{P}(\theta) \le C_2h(\theta).
\]
\end{minipage}

\medskip\smallskip
\begin{minipage}{4in}
(ii) If $\theta$ is a Bernoulli shift then
\[
\sH_p(\theta) = \gamma h(\theta)
\]
where $0 < \gamma < \infty$ is a universal constant.
\end{minipage}
\end{center}

\bigskip
It is not known whether (ii) does not actually hold for all $\theta$.

The definition of $\sH_{P}(\theta)$ easily extends to more general non-singular transformations $\theta$ for which there may be no equivalent invariant probability measure. It is not known whether $\sH_{P}(\theta)$ is a non-trivial invariant for transformations which are not equivalent to transformations with an invariant probability measure. As pointed out by Lewis Bowen to us, the results of \cite{24} may be relevant to this question.

Weaker than the connection to the Kolmorgorov--Sinai entropy, there is also a connection to the Avez entropy of random walks on groups.

Let $\sG$ be a group with a finite generator $g_1,\dots,g_n$ and let $\mu$ be a probability measure with finite support on $\sG$ and let $h(\sG,\mu)$ be the Avez entropy of the random walk on $\sG$ defined by $\mu$. We have the following result \cite{50}:

\bigskip
\begin{center}
\begin{minipage}{4in}
{\em If $h(\sG,\mu) > 0$ then $k^-_{\infty}(\lambda(g_1),\dots,\lambda(g_n)) > 0$ where $\lambda$ is the regular representation of $\sG$ on $\ell^2(\sG)$.}
\end{minipage}
\end{center}

\bigskip
We should also point out to the reader that further results on $k^-_{\infty}$ for Gromov hyperbolic groups $\sG$ and for the entropy of subshifts can be found in the papers \cite{32}, \cite{33}, \cite{34} of Rui Okayasu.

The result about Avez entropy dealt with $k^-_{\infty}(\lambda(g_1),\dots,\lambda(g_n))$. We shall return to this quantity in the next section where we discuss $k_{\sI}(\lambda(g_1),\dots,\lambda(g_n))$ more generally.

\section{Finitely generated groups and supramenability}
\label{sec9}

If $\sG$ is a finitely generated group with generator $\gamma = \{g_1,\dots,g_n\}$ and $\sI$ is a normed ideal, then, which of the following three possibilities takes place
\[
\begin{matrix}
k_{\sI}(\lambda(\gamma)) = 0, \\
0 < k_{\sI}(\lambda(\gamma)) < \infty, \\
k_{\sI}(\lambda(\gamma)) = \infty
\end{matrix}
\]
does not depend on the choice of the generator $\gamma$ and is thus an invariant of the group $\sG$. In particular if $\sI = \sC^-_p$, the number $p_0 \in [1,\infty]$ so that $p \in [1,p_0) \Rightarrow k^-_p(\lambda(\gamma)) = 0$ and $p \in (p_0,\infty] \Rightarrow k^-_p(\lambda(\gamma)) = \infty$ is an invariant of $\sG$, a kind of dimension.

Here are three examples.

\bigskip
\begin{center}
\begin{minipage}{4in}
$1^{\circ}$ \cite{43} If $\sG = {\mathbb Z}^n$ then $k_n^-(\lambda(\gamma)) \in (0,\infty)$.
\end{minipage}

\medskip\smallskip
\begin{minipage}{4in}
$2^{\circ}$ \cite{4} If $\sG$ is the discrete Heisenberg group of $3 \times 3$ upper triangular, unipotent matrices with integer entries, then $k^-_4(\lambda(\gamma)) \in (0,\infty)$.
\end{minipage}

\medskip\smallskip
\begin{minipage}{4in}
$3^{\circ}$ \cite{45} If $\sG$ is a free group on $n \ge 2$ generators, then $k^-_{\infty}(\lambda(\gamma)) \in (0,\infty)$.
\end{minipage}
\end{center}

\bigskip
{\em In view of the special features of $k^-_{\infty}$ and $\sC^-_{\infty}$ it is natural to wonder for which finitely generated groups is $k^-_{\infty}(\lambda(\gamma)) = 0$?}

Here is what we know \cite{54}:

\bigskip
\begin{center}
\begin{minipage}{4in}
(i) {\em if $\sG$ has subexponential growth the $k^-_{\infty}(\lambda(\gamma)) = 0$.}
\end{minipage}

\medskip\smallskip
\begin{minipage}{4in}
(ii) {\em if $k^-_{\infty}(\lambda(\gamma)) = 0$, then $\sG$ is supramenable.}
\end{minipage}
\end{center}

\bigskip
The fact that subexponential growth insures the vanishing of $k^-_{\infty}(\lambda(\gamma))$ is easy. The second assertion uses a recent result of Kellerhals--Monod--R\o rdam \cite{28} which is not easy. For the reader's convenience we include here a few things about the notion of supramenability introduced by Joseph Rosenblatt \cite{38} (we'll stay with finitely generated groups). The group $\sG$ is supramenable if for every subset $\emptyset \neq A \subset \sG$ there is a left invariant, finitely additive measure on the subsets of $\sG$, taking values in $[0,\infty]$, so that $\mu(A) = 1$. In particular supramenable groups are amenable and groups with subexponential growth are supramenable. On the other hand there are amenable groups which are not supramenable. The Kellerhals--Monod-R\o rdam theorem establishes the fact that supramenability of $\sG$ is equivalent to the fact that there is no Lipschitz embedding of a free group on two generators $F_2$ into $\sG$ with respect to the Cayley graph metric. It is also not known whether supramenability and subexponential growth are not actually equivalent properties.

{\it Concerning the class of finitely generated groups for which $k^-_{\infty}(\lambda(\gamma))$ vanishes, it is natural to ask whether it coincides with the class of supramenable groups or with the class of groups with subexponential growth, with the possibility that actually all three classes coincide.} Note also a fourth condition introduced by Monod \cite{31} quite recently and which could be equivalent to some of the preceding.

Amusingly, there is a certain similarity in what we don't know about supramenability and about vanishing of $k^-_{\infty}(\lambda(\gamma))$. We don't know whether the supramenability of $\sG_1$ and $\sG_2$ implies that of $\sG_1 \times \sG_2$. Similarly, we don't know whether vanishing of $k^-_{\infty}$ for $\sG_1$ and $\sG_2$ implies this property for $\sG_1 \times \sG_2$. We should also point out that this is a question specifically for generators of groups, since there are $n$-tuples $\tau$ and $\tau'$ so that $k^-_{\infty}(\tau \otimes I,I \otimes \tau') > 0$, while $k^-_{\infty}(\tau) = k^-_{\infty}(\tau') = 0$ (actually $k^-_p(\tau) = k^-_p(\tau') = 0$ for some given $p > 1$) \cite{54}.

Finally, we should remark that the questions about $k_{\sI}(\lambda(\gamma)) > 0$ are actually questions about functions on $\sG$. If $\ell_{\sI}(\sG)$ denotes the symmetrically normed Banach space on $\sG$, which identifies with the diagonal operators in $\sI$, then $k_{\sI}(\lambda(\gamma)) > 0$ is equivalent to
\[
0 < \inf\left\{\max_{1 \le j \le n} |f(\cdot)-f(g_j\cdot)|_{\sI} \mid f: \sG \to {\mathbb R},\ \mbox{supp } f 
\mbox{ finite, } f(e) = 1\right\}.
\]
The quantity appearing above can be viewed as a generalization of Yamasaki hyperbolicity, which is the special case when $\sI = \sC_p$ \cite{59}.

From this point of view it is easy to see that $k_{\sI}(\lambda(\gamma)) = 0$ is a property of the Cayley graph of $\sG$. Actually, even more:

\bigskip
\begin{center}
\begin{minipage}{4in}
{\em if $\sG_1,\sG_2$ are groups with finite generators $\gamma_1,\gamma_2$ and $\psi: \sG_1 \to \sG_2$ is an injection which is Lipschitz with respect to the Cayley graph metrics, then
\[
k_{\sI}(\lambda(\gamma_2)) = 0 \Rightarrow k_{\sI}(\lambda(\gamma_1)) = 0
\]
(this was used for $\sI = \sC^-_{\infty}$ in \cite{54}).}
\end{minipage}
\end{center}

\bigskip
\section{The commutant $\mbox{\rm mod}$ a normed ideal $\sE(\tau;\sI)$ and its compact ideal $\sK(\tau;\sI)$}
\label{sec10}

To put more structure around the invariant $k_{\sI}(\tau)$, we shall now introduce the commutant mod $\sI$ of $\tau$. We shall assume $\tau = \tau^*$.

\bigskip
\begin{center}
\begin{minipage}{4in}
{\em The commutant {\em mod} $\sI$ of the $n$-tuple of hermitian operators $\tau$ is the subalgebra of $\sB(\sH)$}
\[
\sE(\tau;\sI) = \{X \in \sB(\sH) \mid [X,T_j] \in \sI,\ 1 \le j \le n\}
\]
{\em which is a Banach algebra with isometric involution when endowed with the norm}
\[
\interleave X\interleave = \|X\| + \max_{1 \le j \le n} |[T_j,X]|_{\sI}.
\]
{\em The compact ideal of $\sE(\tau;\sI)$ is}
\[
\sK(\tau;\sI) = \sE(\tau;\sI) \cap \sK
\]
{\em and the Calkin algebra of $\sE(\tau;\sI)$ is the quotient Banach algebra with involution}
\[
\sE/\sK(\tau;\sI) = \sE(\tau;\sI)/\sK(\tau;\sI).
\]
\end{minipage}
\end{center}

\bigskip
Whether $k_{\sI}(\tau)$ vanishes, is finite but $> 0$ or $= \infty$, often appears among the assumptions when studying properties of $\sE(\tau;\sI)$. Actually these 3 situations can be expressed also in terms of approximate units for the compact ideal $\sK(\tau;\sI)$ (\cite{55}, \cite{57}).

\bigskip
\begin{center}
\begin{minipage}{4in}
{\em Assume $\sR$ is dense in $\sI$.}
\begin{itemize}
\item[a)] {\em The following are equivalent}

\medskip\smallskip
\begin{itemize}
\item[(i)] $k_{\sI}(\tau) = 0$

\medskip\smallskip
\item[(ii)] {\em there are} $A_m \in \sK(\tau;\sI)$, $m \in {\mathbb N}$ {\em so that}
\[
\lim_{m \to \infty} \interleave A_mK - K\interleave = \lim_{m \to \infty} \interleave KA_m - K\interleave = 0
\]
{\em for all $K \in \sK(\tau;\sI)$ and} $\interleave A_m\interleave \le 1$, $m \in {\mathbb N}$.

\medskip\smallskip
\item[(iii)] {\em condition (ii) is satisfied and moreover the $A_m$ are finite rank, $0 \le A_m \le I$ and $A_m \uparrow I$ as $m \to \infty$.}
\end{itemize}

\medskip\smallskip
\item[b)] {\em The following are equivalent}
\begin{itemize}

\medskip\smallskip
\item[(i)] $k_{\sI}(\tau) < \infty$

\medskip\smallskip
\item[(ii)] {\em there are $A_m \in \sK(\tau;\sI)$, $m \in {\mathbb N}$ so that}
\[
\lim_{m \to \infty} \interleave A_m K - K\interleave = \lim_{m \to \infty} \interleave KA_m-K\interleave = 0
\]
{\em for all $K \in \sK(\tau;\sI)$ and $\sup_{m \in {\mathbb N}} \interleave A_m\interleave < \infty$.}

\medskip\smallskip
\item[(iii)] {\em condition} (ii) {\em is satisfied and moreover the $A_m$ are finite rank $0 \le A_m \le I$ and $A_m \uparrow I$ as $m \to \infty$.}
\end{itemize}
\end{itemize}
\end{minipage}
\end{center}

\bigskip
In general $\sE/\sK(\tau;\sI)$ is only a Banach-algebra and for purely algebraic reasons it identifies with the $*$-subalgebra $p(\sE(\tau;\sI))$ of the usual Calkin algebra $\sB/\sK$. However, the connection between $k_{\sI}(\tau)$ and approximate unit for $\sK(\tau;\sI)$ has the following somewhat unexpected consequence for $\sE/\sK(\tau;\sI)$.

\bigskip
\begin{center}
\begin{minipage}{4in}
{\em Assuming $\sR$ is dense in $\sI$, we have:}

\medskip\smallskip
\begin{itemize}
\item[a)] {\em if $k_{\sI}(\tau) = 0$ then $p(\sE(\tau;\sI))$ is a $C^*$-subalgebra of $\sB/\sK$ which is isometrically isomorphic to $\sE/\sK(\tau;\sI)$.}

\medskip\smallskip
\item[b)] {\em if $k_{\sI}(\tau) < \infty$ then $p(\sE(\tau;\sI))$ is a $C^*$-subalgebra of $\sB/\sK$, which is isomorphic as a Banach algebra to $\sE/\sK(\tau;\sI)$ (the norms are equivalent).}

\medskip\smallskip
\item[c)] {\em in particular if $\sI = \sC^-_{\infty}$, $p(\sE(\tau;\sC^-_{\infty}))$ is always a $C^*$-algebra canonically isomorphic to the Banach algebra $\sE/\sK(\tau;\sI)$ and if $k^-_{\infty}(\tau) = 0$, the isomorphism is isometric.}
\end{itemize}
\end{minipage}
\end{center}

\bigskip
Concerning the analogy with the usual Calkin algebra, we should also mention the result in \cite{6} about the centre of certain $\sE/\sK(\tau;\sI)$.

\section{Banach space dualities}
\label{sec11}

When $k_{\sI}(\tau) = 0$ and certain additional conditions are satisfied there are many similarities between $(\sK(\tau;\sI),\sE(\tau;\sI))$ and $(\sK,\sB)$ (which can be viewed as the case of $\tau = 0$). For instance, we have the following Banach space duality properties (\cite{52}, \cite{56})

\bigskip
\begin{center}
\begin{minipage}{4in}
$1^{\circ}$ {\em assuming $\sR$ is dense in $\sI$ and in $\sI^d$ and $k_{\sI}(\tau) = 0$, we have that $\sE(\tau;\sI)$ identifies with the bidual of $\sK(\tau;\sI)$}
\end{minipage}

\medskip\smallskip
\begin{minipage}{4in}
$2^{\circ}$ {\em assuming $\sI$ is reflexive and $k_{\sI}(\tau) = 0$, the Banach space $\sE(\tau;\sI)$ has unique predual.}
\end{minipage}
\end{center}

\bigskip
The second property, the uniqueness of predual result, is the analogue of the Grothendieck--Sakai uniqueness (\cite{21}, \cite{39}). It uses a decomposition into singular and ultraweakly continuous part of the functionals on $\sE(\tau;\sI)$ which is analogous to a theorem of M.~Takesaki \cite{41} for von~Neumann algebras. This result then can be used in conjunction with a general result of H.~Pfitzner \cite{36} on unique preduals.

The dual of $\sK(\tau;\sI)$ which is implicit in the above results, under the assumption that $\sR$ is dense in $\sI$ and that $k_{\sI}(\tau) < \infty$ (\cite{52}, \cite{55}) can be identified with $(\sC_1 \times (\sI^d)^n)/\sN$, where $\sN$ is the subspace of elements $(x,(y_j)_{1 \le j \le n})$ so that $x = \sum_{1 \le j \le n} [T_j,y_j]$. The norm on $(\sC_1 \times (\sI^d)^n)$ is
\[
\|(x,(y_j)_{1 \le j \le n}\| = \max\left( |x|_1,\sum_{1 \le j \le n} |y_j|_{\sI^d}\right).
\]
The duality pairings arise by mapping $\sE(\tau;\sI)$ (and hence also $\sK(\tau;\sI)$) into $\sB \times \sI^n$ by
\[
X \to (X,([X,T_j])_{1 \le j \le n}).
\]

\section{Multipliers}
\label{sec12}

If $\sI$ is a normal ideal, we recall that $\sI^{(0)}$ denotes the closure of $\sR$ in $\sI$. Let also
\[
{\tilde \sI} = \left\{X \in \sK \mid \sup_{P \in \sP \cap \sR} |PX|_{\sI} < \infty\right\}
\]
where $\sP$ is the set of hermition projections in $\sB$. The sup in the definition of ${\tilde \sI}$ is the definition of a norm on ${\tilde \sI}$ which is also a normed ideal. The fact that $\sB$ is the multiplier algebra (double centralizer) of $\sK$ has the following analogue in our setting (\cite{52}, \cite{55}):

\bigskip
\begin{center}
\begin{minipage}{4in}
{\em If $\sI$ is a normed ideal and $\tau$ is an $n$-tuple of hermitian operators, so that $k_{\sI}(\tau) < \infty$,, then $\sK(\tau;\sI^{(0)})$ is a closed ideal in $\sE(\tau;{\tilde \sI})$ and $\sE(\tau;{\tilde \sI})$ identifies with the multiplier algebra of $\sK(\tau;\sI^{(0)})$.}
\end{minipage}
\end{center}

\bigskip
\section{Countable degree $-1$-saturation}
\label{sec13}

The countable degree $-1$-saturation property of $C^*$-algebras of I.~Farah and B.~Hart \cite{18}, is a model-theory property which is satisfied by the Calkin algebra and by many other corona algebras. Very roughly, if certain linear relations are satisfied approximately, they can also be satisfied exactly. We have shown in \cite{52} the following:

\bigskip
\begin{center}
\begin{minipage}{4in}
{\em If $\sR$ is dense in $\sI$ and $k_{\sI}(\tau) = 0$, then $\sE/\sK(\tau;\sI)$ has the countable degree $-1$-saturation property of Farah--Hart.}
\end{minipage}
\end{center}

\bigskip
Based on a result of \cite{10}, here is an example of the consequences of degree $-1$-saturation.

\bigskip
\begin{center}
\begin{minipage}{4in}
{\em Assume $\sR$ is dense in $\sI$ and $k_{\sI}(\tau) = 0$. If $\Gamma$ is a countable amenable group and $\rho$ is a bounded homomorphism
\[
\rho: \Gamma \to G L(\sE/\sK(\tau;\sI))
\]
then $\rho$ is unitarizable, that is, there is $s \in G L(\sE/\sK(\tau;\sI))$ so that
\[
s\rho(\Gamma)s^{-1} \subset \sU(\sE/\sK(\tau;\sI)).
\]
}
\end{minipage}
\end{center}

\bigskip
\section{$K$-theory aspects}
\label{sec14}

In case $\sI = \sK$ and $C^*(\tau) \cap \sK = \{0\}$, the algebra $\sE/\sK(\tau;\sI)$ which is the commutant of $p(C^*(\tau))$ in the Calkin algebra $\sB/\sK$, is precisely the Paschke dual of $C^*(\tau)$. The first impulse when encountering the algebra $\sE/\sK(\tau;\sI)$ is to think that they are some kind of smooth subalgebras of the Paschke dual. Already the fact that when $k_{\sI}(\tau) < \infty$, $\sE/\sK(\tau;\sI)$ is a $C^*$-algebra, suggests we are dealing with a quite different situation. Results about the $K$-theory, which show the $K$-theory can be much richer than that of the Paschke dual make this difference quite clear.

To avoid technicalities, we will deal here with generic examples instead of the more general results in the original papers. We shall consider $\tau_n$ the $n$-tuple of multiplication operators by the coordinate functions in $L^2([0,1]^n,d\lambda)$, where $\lambda$ is Lebesgue measure. The $C^*$-algebra of $\tau_n$ is the $C^*$-algebra of continuous functions on $[0,1]^n$. Since $[0,1]^n$ is contractable we infer because of the properties of the Paschke dual construction that the $K_0$-group of $\sE/\sK(\tau_n;\sK)$ is the same as $K_0(\sB/\sK) = 0$.

We shall denote by $\sF_n$ the ordered group of Lebesgue measurable functions $f: [0,1]^n \to {\mathbb Z}$ up to almost everywhere equality and which are in $L^{\infty}([0,1]^n,d\lambda)$, i.e., bounded. Note that $\sF_n$ coincides with the ordered group $K_0((\tau_n)')$ where $(\tau_n)'$ is the von~Neumann algebra which i the commutant of $\tau_n$.

\bigskip
\begin{center}
\begin{minipage}{4in}
{\bf Example 1} \cite{55}. There is an order preserving isomorphism
\[
K_0(\sE(\tau;\sC_1)) \to \sF_1
\]
where for each projection $P$ in $\sM_n(\sE(\tau;\sC_1))$ its $K_0$-class $[P]_0$ is mapped to the multiplicity function of the Lebesgue absolutely continuous part of $P(T_1 \otimes I_n)P$ where $\tau_1 = (T_1)$. Since $k_1(\tau_1) < \infty$ we have that the $K$-group of $\sK(\tau_1;\sC_1)$ and $\sK$ are equal and since $\sE(\tau_1;\sC_1)$ contains a Fredholm operator of index $1$, we have that
\[
K_0(\sE(\tau_1,\sC_1)) \simeq K_0(\sE/\sK(\tau_1,\sC_1)).
\]
\end{minipage}
\end{center}

\bigskip
\begin{center}
\begin{minipage}{4in}
{\bf Example 2} \cite{55}. Assume $\sI \ne \sC_1$ and $\sR$ is dense in $\sI$, then we have
\[
K_0(\sE(\tau_1;\sI)) = 0
\]
and this also implies
\[
K_0(\sE/\sK(\tau_1;\sI)) = 0
\]
since $\sR$ is dense in $\sK(\tau_1,\sI)$ and $K_1(\sK(\tau_1,\sI)) = 0$.
\end{minipage}
\end{center}

\bigskip
\begin{center}
\begin{minipage}{4in}
{\bf Example 3} \cite{55}. If $n \ge 3$ and $\sI = \sC^-_n$, then we have
\[
K_0(\sE(\tau_n,\sC^-_n)) = \sF_n \oplus \sX_n
\]
where the direct summand $\sX_n$ is not known. This can also be stated saying that the map
\[
K_0((\tau_n)') \to K_0(\sE(\tau_n;\sC^-_n))
\]
is an injection and its range is complemented. This uses the results on $\sC^-_n$ perturbations of commuting $n$-tuples of hermitian operators. Using $k^-_n(\tau_n) < \infty$ we also have that $K_0(\sE(\tau_n,\sC^-_n) \simeq K_0(\sE/\sK(\tau_n;\sC^-_n))$.
\end{minipage}
\end{center}

\bigskip
The next example will show that if $n > 1$, there is no analogue of the situation we had when $n = 1$ in Example~2, that is that $K_0$ be trivial if $\sI \supset \sC^-_n$, $\sI \ne \sC^-_n$.

\bigskip
\begin{center}
\begin{minipage}{4in}
{\bf Example 4} \cite{51}. Let $n = 2$ and $\sI = \sC_2$. Then, if $P$ is a projection in $\sM_n(\sE(\tau_2,\sC_2))$, the operator $P((T_1 + iT_2) \otimes I_n)P$ is an operator with trace-class self-commutator and associated with such an operator there is its Pincus principal function $g_{P(T\otimes I_n)P}$ which is in $L^1([0,1]^2,d\lambda)$. The map
\[
[P]_0 \to g_{P(T \otimes I_n)P} \in L^1([0,1]^2,d\lambda)
\]
turns out to be well defined and gives a homomorphism
\[
K_0(\sE(\tau_2,\sC_2)) \to L^1([0,1]^2,d\lambda).
\]
\end{minipage}
\end{center}

\bigskip
One can also show that the range is an uncountable subgroup of $L^1([0,1]^2,d\lambda)$. We refer the reader to (\cite{51}, \cite{46}, \cite{53}) for a discussion about how this homomorphism relates $K_0(\sE(\tau_2;\sC_2))$ to problems on almost normal operators. Note also that the Pincus principal function (\cite{9}, \cite{30}) is related to cyclic cohomology and thus at least some part of $K_0(\sE(\tau_2;\sC_2))$ is related to cyclic cohomology \cite{14}. We should also point out that the algebras $\sE(\tau_2;\sC_2)$ are also related to non-commutative potential theory objects (\cite{11}, \cite{12}, \cite{51}).

\section{The hybrid generalization}
\label{sec15}

In the recent papers \cite{57}, \cite{58} we have shown that the machinery we developed for normed ideal perturbations extends to hybrid perturbations that is $n$-tuples of hermitian operators $\tau$ and $\tau'$ such that $T_j - T'_j \in \sI_j$ where $(\sI_1,\dots,\sI_n)$ is an $n$-tuple of normed ideals. The surprising feature has been that the extension continues to produce sharp results. We shall illustrate this with a few examples of results for commuting $n$-tuples of hermitian operators.

\bigskip
\begin{center}
\begin{minipage}{4in}
{\it If $\tau = (T_j)_{1 \le j \le n}$ is an $n$-tuple of hermitian operators and $(\sI_1,\dots,\sI_n)$ is an $n$-tuple of normed ideals, then $k_{(\sI_1,\dots,\sI_n)}(\tau)$ is defined as the smallest $C \in [0,\infty]$ for which there are $A_m \uparrow I$, $0 \le A_m \le I$ finite rank operators so that}
\[
\lim_{m \to \infty} \max_{1 \le j \le n} |[A_m,T_j]|_{\sI_j} = C.
\]
\end{minipage}
\end{center}

\bigskip
If the $n$-tuple of ideals is $(\sC_{p_1},\dots,\sC_{p_n})$ or $(\sC^-_{p_1},\dots,\sC^-_{p_n})$ we also use the notation $k_{p_1,\dots,p_n}(\tau)$ or $k^-_{p_1,\dots,p_n}(\tau)$ respectively.

\bigskip
\begin{center}
\begin{minipage}{4in}
$1^{\circ}$ \cite{58} {\em Let $\tau$ and $\tau'$ be $n$-tuples of commuting hermitian operators and $p_j \ge 1$, $1 \le j \le n$ so that $p^{-1}_1 + \dots + p_n^{-1} = 1$ and $T_j - T'_j \in \sC^-_{p_j}$ $1 \le j \le n$. Then the absolutely continuous parts $\tau_{ac}$ and $\tau'_{ac}$ are unitarily equivalent.}
\end{minipage}

\medskip\smallskip
\begin{minipage}{4in}
$2^{\circ}$ \cite{57} {\em Let $p_j \ge 1$, $1 \le j \le n$ be so that $p^{-1}_n + \dots + p^{-1}_1 = 1$. Then there is a universal constant $0 < \gamma_{p_1,\dots,p_n} < \infty$ so that if $\tau$ is an $n$-tuple of commuting hermitian operators and $m(x)$, $x \in {\mathbb R}^n$ is the multiplicity function of its Lebesgue absolutely continuous part, we have}
\[
(k^-_{p_1,\dots,p_n}(\tau))^n = \gamma_{p_1,\dots,p_n} \int_{{\mathbb R}^n} m(s)d\lambda(s).
\]
\end{minipage}
\end{center}

\bigskip
\section{Unbounded Fredholm modules}
\label{sec16}

Alain Connes (\cite{13}, see also \cite{15}, \cite{19}) has provided an upper bound for $k^-_n$ based on unbounded Fredholm modules arising in his work on non-commutative geometry.

The unbounded Fredholm module with which one deals here is given by a $*$-algebra of bounded operators on a Hilbert space $\sH$ and an unbounded densely defined self-adjoint operator $D$ so that:

$[D,a]$ when $a \in \sA$, is densely defined and bounded, and $|D|^{-1} \in \sI$, where $\sI$ is a normed ideal.

We refer to $(\sA,D)$ as an unbounded $\sI$-Fredholm-module on $\sH$. (From the early papers \cite{13}, \cite{14} the terminology has been fluid and other related names like spectral triple, $K$-cycle have also been used). Here is the Connes estimate:

\bigskip
\begin{center}
\begin{minipage}{4in}
{\it Let $\tau$ be an $n$-tuple of operators in $\sA$, where $(\sA,D)$ is a $(\sC^-_q)^{dual}$-unbounded Fredholm module, with $q = p(p-1)^{-1}$, then
\[
k^-_p(\tau) \le \beta_p\|[D,\tau]\|(\mbox{Tr}_{\omega}(|D|^{-p}))^{1/p}
\]
where $\beta_p$ is a universal constant and $\mbox{Tr}_{\omega}$ is the Dixmier trace.}
\end{minipage}
\end{center}

\bigskip
Note that the ideal $(\sC^-_q)^{dual}$ is larger than $\sC^-_p$, it is actually $\sC^+_p$ the $(p,\infty)$ ideal on the Lorenz scale. The estimate fits situations involving pseudodifferential operators $D$.

Unless some of the unbounded Fredholm module requirements are relaxed one should not expect a perfect fit between existence of Fredholm modules and $k^-_p$. On the other hand, unbounded Fredholm modules behave well with respect to tensor products, which is not the case for $k^-_p$. For other results around unbounded Fredholm modules and $k_{\sI}$ see also the last part of \cite{45}.

For more on the Connes estimate and non-commutative geometry see (\cite{13}, \cite{15}, \cite{19}).

\section{Sample open problems}
\label{sec17}

Besides the open questions which have come up in our exposition, there are certainly many more. Here are a few we would like to point out.

\bigskip
\begin{center}
\begin{minipage}{4in}
{\bf Problem 1}. Find upper and lower bounds for the universal constants $\gamma_n$, $n \ge 2$ in the formula for $(k^-_n(\tau))^n$ where $\tau$ is an $n$-tuple of commuting hermitian operators. More generally the same question for the universal constants $\gamma_{p_1,\dots,p_n}$ in the hybrid setting is also open. Clearly, it would be of interest to have lower and upper bounds as close to each other.
\end{minipage}
\end{center}

\bigskip
\begin{center}
\begin{minipage}{4in}
{\bf Problem 2}. Does the Farah--Hart degree $-1$-saturation property still hold for $p(\sE(\tau;\sI))$ when the assumption $k_{\sI}(\tau) = 0$ is replaced by $0 < k_{\sI}(\tau) < \infty$? If the answer to the preceding is negative, is there some weaker form of degree $-1$-saturation of $p(\sE(\tau;\sI))$ when $0 < k_{\sI}(\tau) < \infty$? In particular it would be of interest to know the answer to these questions in the case when $\sI = \sC_1$ and $\tau$ is a singleton, a hermitian operator with Lebesgue absolutely continuous spectrum of multiplicity one.
\end{minipage}
\end{center}

\bigskip
\begin{center}
\begin{minipage}{4in}
{\bf Problem 3}. The perturbation entropy of a measure preserving transformation $\sH_{P}(\theta)$ has a natural generalization \cite{47} to an invariant $\sH_{P}(\theta_1,\dots,\theta_n)$ of an $n$-tuple of such transformations
\[
\sup_{\begin{matrix}
\varphi \subset \Phi \\
\varphi \mbox{ finite}
\end{matrix}} k^-_{\infty} (\varphi \cup \{\sU_{\theta_1},\dots,\sU_{\theta_n}\}).
\]
What is the corresponding generalization of the Kolmogorov--Sinai entropy $h(\theta)$ so that
\[
\sH_{P}(\theta_1,\dots,\theta_n)\ {\begin{matrix} {\smile} \\ {\frown} \end{matrix}}\ h(\theta_1,\dots,\theta_n)?
\]
One possible candidate for $h(\theta_1,\dots,\theta_n)$ could be the supremum over finite partitions $\sP$ of the probability measure space of
\[
\liminf_{m \to \infty} m^{-1}H(\sP_m)
\]
where $\sP_n$ is defined recursively by $\sP_1 = \sP$ and $\sP_{m+1} = \sP_m \vee \theta_1\sP_m \vee \dots \vee \theta_n\sP_m$.
\end{minipage}
\end{center}

\end{document}